\documentclass[12pt]{article}
\usepackage{combelow}

\newcommand{\Kn}[4]{        
	\begin{tikzpicture}[rotate=#2, scale=3]
	\foreach \i in {1,...,#1}
	{
		\path (360/#1*\i:1cm)  node (X\i){$\i$};
		#4 (X\i) circle (#3pt);
	}
	\foreach \i in {1,...,#1}
	{
		\foreach \j in {1,...,#1}
		{ 
	 	\pgfmathparse{int{mod(\i+\j,3)}}
	     \ifnum\pgfmathresult=0{\draw [color=blue](X\i)--(X\j);}
			\else{
				\ifnum\pgfmathresult=1{\draw[color=red](X\i)--(X\j);}
				\else{\draw [color=green](X\i)--(X\j);}
				\fi
			}
			\fi
		}
	}
	\end{tikzpicture}
}

\usepackage[bookmarksnumbered,colorlinks,plainpages,backref=false]{hyperref}

\usepackage{tikz,pgf}
\usepackage{latexsym,amsmath,graphicx}
\usepackage[latin2]{inputenc}
\usepackage[T1]{fontenc}
\usepackage{t1enc}
\usepackage[magyar]{babel}

\def\idez#1{{\ideze}#1''}
\def\ideze{\setbox0=\hbox{\lower1.38ex\hbox{''}}\dp0=0pt\box0}

\newenvironment{proof}{\bizo}{\kocka}

\def\bizo{\par\smallskip\noindent{\textbf{Bizonyítás.\ }}}
\def\kocka{\nobreak\hspace*{\fill}\nobreak{\piros{$\Box$}}\par\smallskip}

\newcommand{\keret}[1]{\bigskip\noindent\fbox{\parbox{0.98\textwidth}{#1}}}

\newcommand{\kockaa}{\nobreak\hspace*{\fill}\nobreak$\Box$\par\smallskip}

\newcommand{\bizof}[1]{\par\noindent{\textbf{A #1. tétel bizonyítása\ }}}
\newcommand{\no}{\medskip\noindent}

\newcommand{\kek}[1]{\textcolor{blue}{\textbf{#1}}}
\newcommand{\piros}[1]{\textcolor{red}{\textbf{#1}}}

\begin{document}

\begin{center}\kek{\textsc{\LARGE Rend a rendezetlenségben -- gráfelméleti példákkal}}

\bigskip\medskip
\kek{\textsc{\LARGE Order in the chaos with examples from graph theory}}
\end{center}

\medskip
\textit{In all chaos there is a cosmos, in all disorder a secret order.} (C. G. Jung)

\bigskip\medskip\bigskip\medskip \centerline{\textbf{\Large English abstract}}

\bigskip\medskip \no In randomly created structures (be they natural or artificial) very often there exist ordered substructures. In this Hungarian language scientific essay we will present some of such structures in graph theory.

\no\textbf{Rédei's theorem:} In a digraph obtained by arbitrarily directing the edges of a complete graph with at least two vertices, there is always a directed Hamiltonian path.

\no \textbf{Ramsey's theorem:} Let $ R (m, k) $ be the smallest natural number for which in every graph with $n$ vertices, where $ n \ge R (m, k) $ there is a complete subgraph $ K_m $  or in its complement a complete subgraph $ K_k $. $ R (m, k) $ always esists for all positive  integer $m$ and $n$.

\no The \textit{Erdős--Szekeres theorem} states that:
$R(m,k)$ always exists for every  $m,k\in \mathbf{N^*}$ and 
\[R(m,k)\le \binom{m+k-2}{m-1}\]
gives an upper bound for it. Exact values for $R(m,k)$ are known only for small values of $m$ and $n$.-

\no\textbf{Turán's theorem:}
Any simple graph with $n$ vertices that does not contain a complete subgraph with  $(k + 1)$ vertices may contain  at most
\[e \le \frac {1} {2} \left (n ^ 2-r ^ 2 \right) \frac {k-1} {k} + \frac {r (r-1)} {2} \]
edges,  where $ n = hk + r $, $ 0 \le r <k. $

Extreme graph (for which there is an equality in the above formula) is $ K_ {n_1, n_2, \ldots, n_k} $, where there are $r$ independent sets with $ h + 1 $ vertices and $k-r$ independent  sets with  $ h $ vertices (So, $ r $ of $n_1, n_2, \ldots, n_k $ equals $ h + 1 $ and $ k-r $ of them equals $ h $). 

\no \textbf{De Bruijn  graphs:}
De Bruijn graphs are very regular and harmonious, but they have a lot of directed Hamiltonian cycles, and these  are not necessarily in some order. The question is whether there may be a certain number of Hamiltonian cycles among them that have no common arcs.

\no\textit{Conjecture of Johny Bond and Antal Iványi:}
If $n \ge 2$ and $m \ge 2$, then in the De Bruijn graph $B(n, m)$ there exists $n-1$ arc-independent directed Hamiltonian cycles. 

\no We also formulate a perhaps easier-to-prove version of this conjecture. 
 
\newpage

\begin{center}
	\kek{\textbf{\LARGE Rend a rendezetlenségben -- gráfelméleti példákkal}}

\medskip\kek{Kása Zoltán}
\end{center}

\bigskip\bigskip

\begin{minipage}{5cm}
\bigskip
$\ $
\end{minipage}	
\begin{minipage}{8cm}
\textit{Minden káoszban van kozmosz (ékesség, rend, tisztesség) is, minden rendetlenségben rend is rejlik, s minden önkényben folyton ott a törvény, mert mindaz, ami hat, ellentéten alapszik.} (C. G. Jung)
\end{minipage}	

\bigskip\bigskip
\noindent
{\small 
Ez az írás\footnote{A Magyar Tudomány Ünnepe konferencián Kolozsvárott 2021. november 11-én elhangzott előadás szerkesztett változata} egyféle gráfelméleti esszészerűség\footnote{Samuel Johnson szerint az esszé: ,,az agy gondolatainak szabad folyása; egyetlen, nyers mű, s nem szabályszerű, rendezett alkotás.''
	}. Azt vizsgáljuk, hogy egy struktúra, legyen az bármennyire is véletlenszerűen, rendezetlenül létrehozva, tartalmazhat-e valamilyen rendezett részstruktúrát. Aki nem járatos a gráfelméletben, a függelékben megtalálja az alapvető fogalmakat. 
}

\bigskip\bigskip
\noindent Ha alaposan megfigyelünk bármilyen rendezetlen struktúrát, észrevehetjük, hogy valahol, valamilyen módon van benne rend. Ezt a rendet keressük a gráfelméletben!

Szemerédi Endre\footnote{Szemerédi Endre (1940. augusztus 21. --) Abel- és Széchenyi-díjas magyar matematikus, egyetemi tanár, a Magyar Tudományos Akadémia rendes tagja. Nemzetközi tudományos ismertségre kombinatorikai, számelméleti és algoritmuselméleti kutatásaival tett szert.} nyilatkozta egy interjúban: \textit{,,A véges objektumokban mintázatokat, különböző alakzatokat keresünk, illetve azt, hogy milyen feltételek mellett jönnek létre bizonyos alakzatok -- ez az egyik alapkérdés. Egy kicsit nagyképű vagy filozofikus megfogalmazás szerint azt szeretnénk bizonyítani, hogy minden káoszban van rend. Tehát ha ön direkt rosszindulatúan ad nekem egy struktúrát, abban is meg lehet találni a rendnek tekinthető részleteket. Legyen mondjuk hat pont, és én azt mondom, hogy ezeket ön kösse össze kék és piros ceruzával összevissza, létrehozva egy élekből és csomópontokból álló, úgynevezett teljes gráfot. Mi lenne ebben egy szép alakzat? Mondjuk három pont, amelyeket összekötő élek azonos színűek. Igazolni lehet, hogy az ördög -- vagy a rosszindulatú újságíró -- bárhogyan is színezi ki a gráfot, könnyen lehet találni benne három ilyen pontot.''}

\newpage
\no\centerline{\piros{Rédei tétele}}

\no Egy utat \emph{Hamilton-útnak} nevezünk, ha tartalmazza a gráf összes csúcsát. Egy kör, amelyik tartalmazza a gráf minden csúcsát, \emph{Hamilton-kör}. 
Irányított gráfokban \textit{irányított Hamilton-útról} beszélünk. R\'edei László\footnote{Rédei László (1900--1980)  Kossuth-díjas matematikus, szegedi  egyetemi tanár, a magyar absztrakt algebrai iskola megalapozója. A Magyar Tudományos Akadémia rendes tagja.} következő érdekes eredménye szerint egy véletlenszerűen  felépített teljes gráfban létezik bizonyos rendezett alstruktúra.

\keret{\textbf{1. tétel.} 

\textit{\kek{Egy legalább két csúcsú teljes gráf éleinek tetszőleges irányításával kapott gráfban mindig van irányított Hamilton-út.}}}

\medskip\begin{proof}
	Tekintsünk a teljes gráfban az élek irányítása után egy leghosszabb $L$ irányított utat (kezd\H ocsúcsa $a$, végcsúcsa $b$):
	
	\bigskip\bigskip
	\qquad\quad\begin{tikzpicture}[scale=1.5]
	\tikzstyle{every node}=[draw,shape=circle];
	\path (0,0) node (v1) {$a$};
	\path (0.75, 0) node (v2) {$\ \ $};
	\path (1.5, 0) node (v3) {$\ \ $};
	\path (2.25, 0) node (v4) {$\ \ $};
	\path (3,0) node (v5) {$p$};
	\path (3.75,0) node (v6) {$q$};
	\path (4.5, 0) node (v7) {$\ \ $};
	\path (5.25, 0) node (v8) {$\ \ $};
	\path (6,0) node (v9) {$b$};  
	\path (3,-1.5) node (v10) {$r$}; 
	
	\draw [->,-latex] (v1)--(v2);
	\draw [->,-latex] (v2)--(v3);
	\draw [->,-latex] (v3)--(v4);
	\draw [->,-latex] (v4)--(v5);
	\draw [->,-latex] (v5)--(v6);
	\draw [->,-latex] (v6)--(v7);
	\draw [->,-latex] (v7)--(v8);
	\draw [->,-latex] (v8)--(v9);
	\draw [->,-latex,color=red] (v1)--(v10);
	\draw [->,-latex,color=red] (v5)--(v10);
	\draw [->,-latex,color=blue] (v10)--(v6);
	\draw [->,-latex,color=blue] (v10)--(v9);
	
	\tikzstyle{every node}=[];
	\draw (7,0) -- (7,0) node[anchor=east] {{\LARGE $L$}};
	\end{tikzpicture}

	\no Ha ez nem Hamilton-út, akkor létezik egy $r$ csúcs a gráfban, amelyik nincs rajta ezen az uton. 
	Ekkor létezik az $(a,r)$  irányított él (különben, ha az irányítás  $(r,a)$ lenne, akkor $L$ nem a leghosszabb irányított út), hasonlóképen létezik az $(r,b)$ él is. Az $L$ út minden csúcsa éllel van összekötve $r$-rel, mivel az eredeti gráf teljes gráf. Mivel nem lehet  $L$ minden bels\H o csúcsa $r$ felé irányított éllel öszzekötve $r$-rel, léteznie kell a $(p,r)$ és $(r,q)$ irányított éleknek, de ekkor  az $a,\ldots p,r,q, \ldots, b$ irányított út hosszabb mint $L$, ami ellentmondás. Tehát $L$ Hamilton-út.
\end{proof}

\no Ebben a tételben az a szép, hogy a teljes gráfot tetszőleges módon irányítva, egy látszólag teljesen rendezetlen gráfot kapunk, de az így kapott  gráfban mégis  mindig van irányított Hamilton-út.

\bigskip
Az alábbi két példában az élek egy-egy irányítását láthatjuk, és pirossal a Hamiltom-utakat. Itt látható, hogy az élek tetszőleges irányításával irányított Hamilton-kör nem mindig kapható. Mindkét esetben a 3-as csúcsba  csak befutó élek vannak, tehát Hamilton-kör léte kizárható.  

\medskip\qquad\qquad\begin{tikzpicture}[scale=1.5]
\tikzstyle{every node}=[draw,shape=circle];
\path (0,0) node (v1) {1};
\path (0,2) node (v2) {2};
\path (2,2) node (v3) {3};
\path (2, 0) node (v4) {4};

\draw [->, -latex, color=red] (v2) -- (v1);
\draw [->, -latex] (v2) -- (v3);
\draw [->, -latex] (v2) -- (v4); 
\draw [->, -latex, color=red] (v1) -- (v4); 
\draw [->, -latex] (v1) -- (v3); 
\draw [->, -latex, color=red] (v4) -- (v3); 
\end{tikzpicture}\qquad \qquad	
\begin{tikzpicture}[scale=1.5]
\tikzstyle{every node}=[draw,shape=circle];
\path (0,0) node (v1) {1};
\path (0,2) node (v2) {2};
\path (2,2) node (v3) {3};
\path (2, 0) node (v4) {4};

\draw [->, -latex] (v2) -- (v1);
\draw [->, -latex] (v2) -- (v3);
\draw [->, -latex, color=red] (v2) -- (v4); 
\draw [->, -latex, color=red] (v4) -- (v1); 
\draw [->, -latex, color=red] (v1) -- (v3); 
\draw [->, -latex] (v4) -- (v3); 
\end{tikzpicture} 

\bigskip
\no\centerline{\piros{Ramsey-típusú széls\H oérték-feladatok}}

\bigskip\noindent Bizonyítsuk be, hogy egy legalább hattagú  
társaságban mindig vannak hárman, akik kölcsönösen ismerik egymást vagy kölcsönösen nem ismerik egymást (az ismeretséget kölcsönösnek tekintjük). Legyenek a gráf csúcsai a társaság emberei.  Jelöljük piros éllel, ha két ember ismeri egymást, és kék éllel ha nem. Válasszuk ki az $x_1$ csúcsot és az $x_2$, $x_3$, $x_4$, $x_5$, $x_6$ csúcsokat. Mivel $x_1$ öt éllel van összekötve a többi csúccsal, és ezek két színnel vannak kiszínezve, az öt él közül, legalább három azonos színű, például piros. Az ábrán feketével szaggatottan rajzolt él mindegyike  lehet akár piros, akár kék.

\no\centerline{\begin{tikzpicture}[scale=1.5]
\tikzstyle{every node}=[draw,shape=circle];
\path (3,1) node (v1) {$x_1$};
\path (1,0) node (v2) {$x_2$};
\path (2,0) node (v3) {$x_3$};
\path (3,0) node (v4) {$x_4$};
\path (4,0) node (v5) {$x_5$};
\path (5,0) node (v6) {$x_6$};
\draw  [color=red,thick]
(v1) -- (v2)
(v1) -- (v3)     
(v1) -- (v4);
\draw   [color=blue,thick]
(v1) -- (v5)
(v1) -- (v6); 
\draw [dashed,thick] 
(v2) -- (v3)
(v3) -- (v4)
(v2) .. controls (2,-0.5) .. (v4);    
\end{tikzpicture}}

\no Ha az $\{x_2,x_3 \}$ él piros, akkor 
a gráfban van piros háromszög ($x_1, x_2,x_3$). Ugyanez igaz az $\{x_3,x_4 \}$ és $\{x_2,x_4 \}$ élekre is. Ha viszont ezek egyike sem piros, akkor mindhárom kék (az ábrán szaggatottan), és akkor van a gráfban kék háromszög. 

\no A feladat úgy is megfogalmazható, hogy bármely legalább hat csúcsú gráfban van a komplementerében van háromszög. (A piros élek a gráfban vannak, a kékek a komplementerében.)

\no A nagyon fiatalon elhunyt Frank Ramsey\footnote{Frank P. Ramsey (1903--1930) brit matematikus, filozófus, közgazdász} matematikus írta, hogy:

\no \kek{\idez{egy elég nagy rendszerben, mégha rendezetlen is, kell lenni valamilyen rendnek} }  

\no Az előbbi feladat (egy legalább hattagú  
társaságban mindig vannak hárman, akik kölcsönösen ismerik egymást vagy kölcsönösen nem ismerik egymást) általánosítása Ramsey nevéhez fűz\H odik.
 
\no Jelöljük $K_n$-nel az $n$-csúcsú teljes gráfot.

 \no 
 \kek{\emph{Legyen $R(m,k)$ az a  legkisebb természetes szám, amely\-re egy  $n\ge  R(m,k)$ csúccsal rendelkez\H o gráfban van  $K_m$ részgráf vagy a komplementerében van $K_k$ részgráf. 
 	Ext\-rém gráfoknak nevezzük azokat   az  $R(m,k)-1$ csúcsú gráfokat, amelyekre ez a tulaj\-donság nem igaz.}
 }

\medskip\no 
Például: $R(3,3)=6$. Extrémgráf az 5 csúcsú kör. Az alábbi példákban a kék élek a gráfhoz tartoznak, a pirosak  a komplementer gráfhoz. 

\bigskip
\qquad\qquad\begin{tikzpicture}[scale=1.5]
\pgfmathsetmacro{\n}{6};       
\foreach \i in {1,...,\n}
{  
	\path (360/\n*\i:1cm) node (X\i) { }; 
	\draw (X\i) circle  (2pt);               
}
\draw [color=blue] (X1) -- (X2);	
\draw [color=red] (X1) -- (X3);
\draw [color=red] (X1) -- (X4);
\draw [color=blue] (X1) -- (X5);
\draw [color=red] (X1) -- (X6);

\draw [color=blue] (X2) -- (X3);
\draw [color=blue] (X2) -- (X4);
\draw [color=red] (X2) -- (X5);
\draw [color=blue] (X2) -- (X6);

\draw [color=blue] (X3) -- (X4);
\draw [color=blue] (X3) -- (X5);
\draw [color=red] (X3) -- (X6);

\draw [color=blue]  (X4) -- (X5);
\draw [color=red]  (X4) -- (X6);

\draw [color=blue]  (X5) -- (X6);
\end{tikzpicture}
\qquad\qquad\quad
\begin{tikzpicture}[rotate=18,scale=1.5]
\pgfmathsetmacro{\n}{5};       
\foreach \i in {1,...,\n}
{  
	\path (360/\n*\i:1cm) node (X\i) { }; 
	\draw (X\i) circle (2pt);               
}

\draw [color=blue] (X1) -- (X2);
\draw [color=red] (X1) -- (X3);
\draw [color=red] (X1) -- (X4);

\draw [color=blue] (X2) -- (X3);
\draw [color=red] (X2) -- (X4);
\draw [color=red] (X2) -- (X5);

\draw [color=blue] (X3) -- (X4);
\draw [color=red] (X3) -- (X5);

\draw [color=blue] (X4) -- (X5);
\draw [color=blue] (X5) -- (X1);	
\end{tikzpicture}

\qquad \ $K_6$ kék háromszöggel \qquad  extrémgráf (nincs egyszín\H u háromszög)

\bigskip
\no Könny\H u bizonyítani, hogy:

$R(1,k)=R(k,1)=1$

$R(2,k)=R(k,2)=k$   

\no Általánosan: $R(m,k)=R(k,m)$.

\no Be lehet bizonyítani, hogy $R(m,k)$ létezik, minden nullától különböző $m$ és $n$ természetes számra.

\keret{\textbf{2. tétel.}
\textit{\no\kek{\textbf{a)} Ha $R(m-1, k)$ és $R(m,k-1)$ létezik, akkor $R(m,k)$ is létezik, és}\\		
		\centerline{\kek{$ R(m,k)\le R(m-1,k)+R(m,k-1) .$}}\\
		\no\kek{\textbf{b)} Ha $R(m-1,k)=2p$ és $R(m,k-1)=2q$,  $p,g\in \mathbf{N^*}$, akkor} 
\\		
			\centerline{\kek{$ R(m,k)< R(m-1,k)+R(m,k-1).$}}
	}
\no \kek{ahol $\mathbf{N^*}$ a pozitív természetes számok halmaza.}
}

\medskip\no Innen következik az Erdős\footnote{Erdős Pál (1913--1996) világhírű magyar matematikus, a Magyar Tudományos Akadémia tagja, aki főleg számelmélettel, kombinatorikával és gráfelmélettel foglalkozott.}--Szekeres\footnote{Szekeres György (1911--2005) magyar-ausztrál matematikus, a Magyar Tudományos Akadémia tagja, aki  gráfelmélettel, kombinatorikával és számelmélettel foglalkozott.}-tétel:

\keret{\textbf{3. tétel}
 
	\textit{\no \kek{$R(m,k)$ létezik tetsz\H oleges  $m,k\in \mathbf{N^*}$ értékekre, és 
		\[R(m,k)\le \binom{m+k-2}{m-1}.\]}
}
}

\medskip
\no \kek{Néhány ismert érték:}

\noindent\url{http://mathworld.wolfram.com/RamseyNumber.html}

\bigskip
\begin{center}
{\footnotesize
	\no
	\begin{tabular}{|c||c|c|c|c|c|c|c|c|}\hline
		$m \backslash k$ & 3 & 4 & 5 & 6 & 7 & 8 & 9 & 10 \\ \hline\hline
		3 & 6 &9&14&18&23&28&36&[40,43]\\ \hline 
		4 &&18&25&[35,41]&[49,61]&[56,84]&[73,115]&[92,149]\\ \hline
		5 &&&[43,49]&[58,87]&[80,143]&[101,216]&[125,316]&[143,442]\\ \hline
		6 & & & & [102,165]&[113,298]&[127,495]&[169,780]&[179,1171]\\ \hline
		7 &&&&&[205,540]&[216,1031]&[233,1713]&[232,2826]\\ \hline
		8 &&&&&&[282,1870]&[317,3583]&[377,6090]\\ \hline
		9 &&&&&&&[565,6588]&[580,12677]\\ \hline
		10&&&&&&&&[798,23556]\\ \hline
	\end{tabular}
}
\end{center}
	
\no Erdős Pál bizonyította, hogy
\[R(k,k)> 2^{k/2} \]

\no
\kek{Megoldatlan:} \[\textrm{Létezik } \lim_{k\rightarrow \infty}R(k,k)^{1/k}?\]

\no
\kek{Általánosítás: } Legyen $R(k_1, k_2, \ldots, k_n)$ az a legkisebb természetes szám, amelyre igaz, hogy ha egy legalább ennyi csúcsból álló teljes gráf éleit $n$ színnel kiszínezzük, mindig létezik egy olyan $i$ szín  ($1\le i\le n$), hogy a gráfban van $k_i$ csúcsú teljes részgráf, amelynek élei az $i$ színnel vannak kiszínezve. 
\[R(k_1+1, k_2+1, \ldots, k_n+1) \le \frac{(k_1+k_2+\ldots +k_n)!}{k_1!k_2!\cdot \ldots \cdot k_n!}\]
Erd\H os Pál bebizonyította, hogy:
\[R(\underbrace{3,3,\ldots, 3}_{r\textrm{-szer} })\; \le\; r!\sum_{k=0}^{r}{\frac{1}{k!}}+1\]

\no Ugyancsak az Erdős Pál sejtése (amely nincs bizonyítva még), hogy mindig igaz az egyenl\H oség. 

\no Ismert értékek: $R(3,3)=6$, $R(3,3,3)=17$.

\no A következő ábrán a $K_{17}$ gráf éleit három színnel színeztük (piros, kék és zöld). Mivel a színezés nem véletlenszerűen történt (az $i$ és $j$ csúcsokat kék éllel kötöttük össze, ha összegük 3-mal való osztási maradéka 0, piros, ha a maradék 1, és zöld, ha a maradék 2), minden színből van háromszög (pl. kék: 3, 9, 15; piros: 5, 11, 17, zöld: 4, 10, 16). 

\begin{center}
\Kn{17}{0}{2}{\draw} 	
\end{center}

\keret{
\textbf{4. tétel. }

	{\kek{$R(3,4)=9.$}}
}

\medskip\no\textbf{Bizonyítás. }
Az 2. tételb\H ol: 

\medskip
$R(3,4)\le R(2,4)+R(3,3)=4+6$, 

\no és mivel mindkettő páros: $R(3,4)< 10$. 

\no Be kell bizonyítani, hogy $R(3,4)> 8.$ 

\no Legyen $H_{3k-1}$ a következ\H o gráf: csúcsai $a_1,a_2, \ldots, a_{3k-1}$ egy körön szabályosan elhelyezve, élekkel kötjük össze azokat, amelyeknek távolsága nagyobb, mint a körbe írható sza\-bályos háromszög oldala. Így pl. $a_1$ össze van kötve $a_{k+1}, \ldots, a_{2k}$ csúcsokkal. 

\medskip\no
\qquad\qquad\qquad\qquad\qquad\includegraphics[scale=0.5]{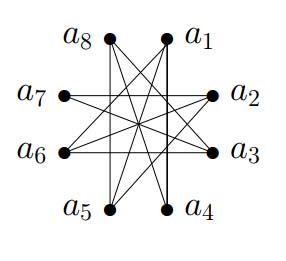}

\qquad \qquad \qquad  \qquad\qquad\qquad\qquad  $H_8$

\bigskip\no  $H_{3k-1}$ nem tartalmaz háromszöget, és legfeljebb $k$ független csúcsa van (pl. $a_1, a_2, \ldots a_k$), ezért 

\medskip $R(3,k+1)\ge 3k$. 

\bigskip\no Ha $k=3$, akkor $R(3,4)\ge 9$. Ezért  $R(3,4)=9.$

\no Tehát \kek{$ 3(k-1) \le R(3,k)\le \frac{k(k+1)}{2}.$}

\no A $H_{3k-1}$ gráfokat \textit{Andrásfai-gráfoknak} nevezzük, mivel Andrásfai Béla\footnote{Andrásfai Béla (1931. február 8. --) matematikus, a budapesti Műszaki Egyetem nyugdíjas docense, az Andrásfai-gráf névadója} írta le először.

\newpage
\no Andrásfai-gráfok $n=1, 2,  \ldots ,9$ értékekre:

\includegraphics[scale=0.5]{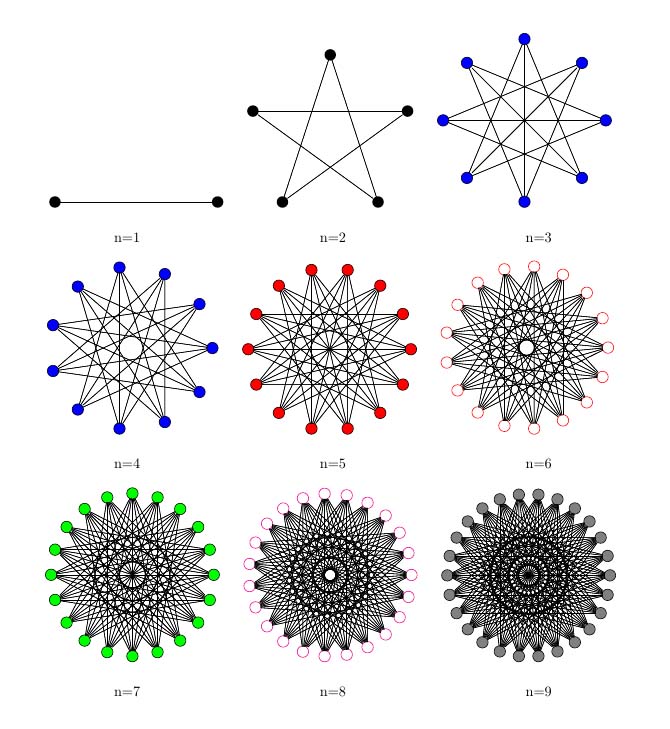}

\bigskip\bigskip
\no\centerline{\piros{Turán-típusú széls\H oérték feladatok}}

\no Ebben az esetben azt vizsgáljuk meg, hogy ha egy gráfban nincs bizonyos nagyságú teljes gráf, akkor maximálisan hány éle lehet a gráfnak.
Itt extrémgráf az a gráf lesz, amely az adott feltételek mellett a legtöbb élt tartalmazza. A következő eredmény Turán Pál\footnote{Turán Pál (1910--1976) matematikus, a Magyar Tudományos Akadémia tagja, a számelmélet, a gráfelmélet és a klasszikus analízis területén ért el jelentős eredményeket.} nevéhez fűződik. 

\keret{
	\textbf{5. tétel.}

\emph{\kek{Bármely  $n$ csúcsú egyszer\H u gráf, amely nem tartalmaz $(k+1)$ csúcsú teljes részgráfot, legfeljebb
	\[e\le \frac{1}{2}\left(n^2-r^2 \right)\frac{k-1}{k}+\frac{r(r-1)}{2}\]
	élt tartalmazhat, ahol  $n=hk+r$, $0\le r<k.$}\\
	\kek{Extrémgráf (amelyre fennáll az egyenl\H oség):  $K_{n_1,n_2, \ldots , n_k}$, ahol  $r$ független csúcshalmaz $h+1$ csúcsból áll és  $k-r$ független csúcshalmaz $h$ csúcsból áll (Tehát, $n_1, n_2, \ldots , n_k$ közül  $r$ egyenl\H o  $(h+1)$-gyel és $k-r$ egyenl\H o $h$-val).}
 }
}

\medskip\no
\textbf{Példa.}

\no Vizsgájunk meg az élek maximális számát egy 5 csúcsú gráfban, amelyben nincs háromszög. Ekkor a tétel szerint $k=2$, és mivel $5= 2\cdot 2 + 1$, ezért $h=2, r=1$. Az extrémgráf, amely nem tartalmaz háromszöget és maximélis számú éle van, a következő (két rajzolatban):

\qquad\qquad\includegraphics[scale=0.75]{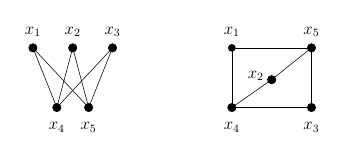}

\no Egy másik példa a következ\H o, ahol  a csúcsok száma 7, és a  gráfban nincs $K_4$.  
 Ekkor $n=7$, $k=3$, és mivel  $7=2\cdot 3+1$, akkor  $h=2$ és $r=1$. 

\no $n_1=3, n_2=2, n_3=2$, és az extrémgráf $K_{3,2,2}$ (amelyik izomorf $K_{2,3,2}$-vel). Az extrémgráf két rajzolatban:

\bigskip
\qquad\qquad\qquad\qquad\includegraphics[scale=0.25]{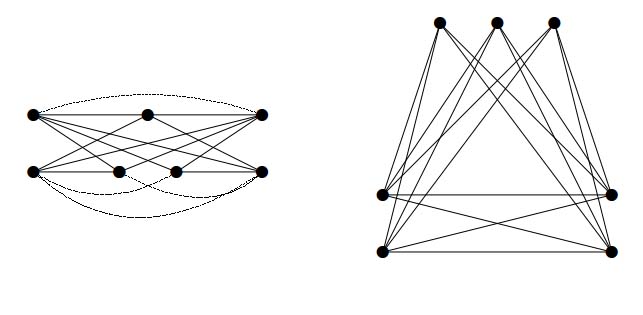}

\no\textbf{Példa.}	
$n=13,\  k=4$,\ $13=\piros{3}\cdot 4+\piros{1}$,\  tehát  $h=3,\   r=1$. Az extrémgráf $K_{4,3,3,3}$:

\begin{center}
	\includegraphics[scale=0.8]{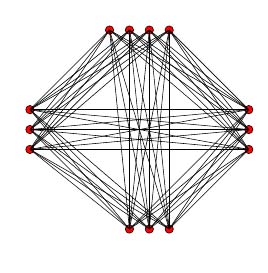}
\end{center}

\no A Szemerédi-féle regularitási lemma olyan gráfelméleti  tétel\footnote{A tétel nevében megmaradt a lemma szó, mivel Szemerédi Endrének annak idején lemmaként volt rá szüksége.}, amely szintén a rendet keresi a rendezetlenségben, azaz igen nagy véletlenszerű gráfokban talál valamilyen rendezett struktúrát. A lemma szerint minden, kellően nagy gráf csúcsainak a halmaza felosztható olyan hasonló méretű részhalmazokra, amelyek esetében a részhalmazok közötti élek csaknem véletlenszerűek.
 
\bigskip\centerline{\piros{De Bruijn-gráfok}}

\no Azt, hogy a rendezetlen nagyméretű struktúrában van rendezett alstruktúra, azt láttuk. Az a kérdés, hogy egy rendezett struktúrában lehet-e rendezetlen alstruktúra. Egészen kisméretú biztos lehet, de nagyméretű lehet-e?

\no A következőkben a De Bruijn-gráfokról lesz szó, amelyek nem a rendezetlenség, hanem éppen ellenkezőleg a rendezettség, a rend, a harmónia képviselői. 

\newpage
\no\kek{Szógráf} alatt olyan gráfot értünk, amelyben a csúcsok $m$ hosszúságú szók\footnote{Itt \textit{szó} alatt egy adott ábécé (jelek véges, nem üres halmaza) egymás mellé helyezett betűiből képzett füzért értjük. Például, ha az ábécé \{$a, b, 1, 2, 3$\}, akkor szók lehetnek: $ab$, $aa$, $1ba2$, $221$ stb.}, irányított él van az 
$a_1a_2\ldots a_m$ csúcsból a  $b_1b_2\ldots b_m$ csúcsba, ha

\centerline{$a_2\ldots a_m=$   $b_1b_2\ldots b_{m-1}$.} 

\medskip
\no Ha $V=A^m$ (azaz, a csúcsok halmaza az összes $m$ hosszúságú szó), ahol $A=\{a_1,a_2, \ldots a_n\}$, akkor a gráf neve \kek{De Bruijn-gráf,} jelölése pedig \kek{$B(n,m)$.} Például, ha $A=\{0, 1\}$, akkor $n=2$, és ha $m=3$, akkor a következő De Bruijn-gráfot kapjuk:

\begin{center}
\includegraphics[scale=0.06]{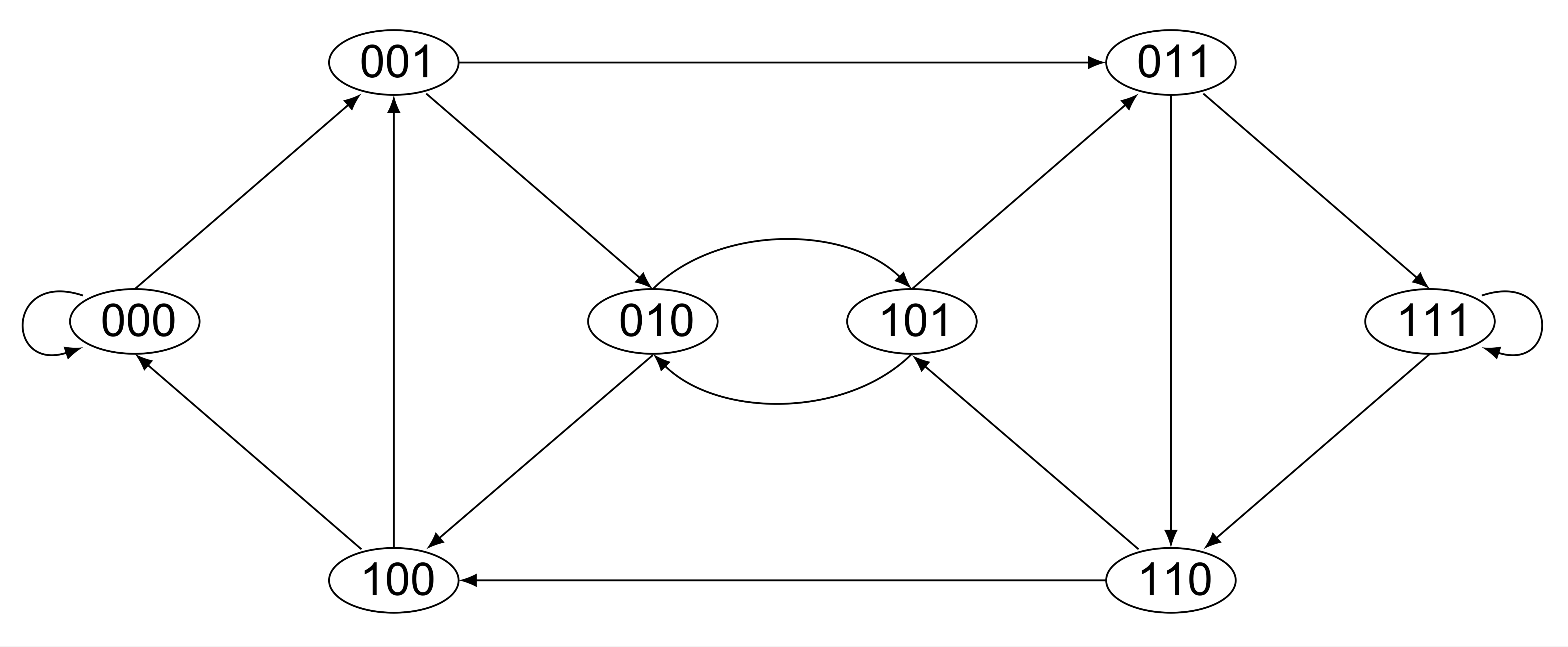}

B(2,3)
\end{center}

\no Ebben a gráfban két irányított Hamilton-kört is találunk. Ezek:

\kek{$000, 001, 011, 111, 110, 101, 010, 100, 000$} és

\kek{$000, 001, 010, 101, 011, 111, 110,  100, 000$}.   

\no Ha elhagyjuk a vesszőket, az egymás melletti szavakat összecsúsztatjuk (a közös részeket csak egyszer írjuk le), és az utolsó betűt is elhagyjuk, akkor ezeket kapjuk: 
\kek{$0001110100$}\  és\  \kek{$0001011100$}, és az így kapott szókat De Bruijn-szóknak nevezzük, és ezek tulajdonképpen De Bruijn-utaknak felelnek meg, de, ha megegyezünk abban, hogy az utolsó betű után mindig $0$ következik, akkor tekinthetjük úgy, hogy De Bruijn-körökről van szó.

\no A $B(n,m)$ De Bruijn-gráfban összesen
$\displaystyle\frac{(n!)^{n^{m-1}}}{n^m}$	
irányított Hamilton-kör van. Ez a szám, ha $n$ és $m$ nő, igen gyorsan növekszik.\\

\begin{center}
$\begin{array}{|l|r|}\hline
\textrm{\ \ gráf}& \textrm{Hamilton-körök száma} \\ \hline
   B(2,3)& 2\\ \hline
   B(3,2)& 24\\ \hline
   B(3,3)&   373248 \\ \hline
   B(3,4)&   13824\cdot 10077696^{3} \\ \hline
\end{array}$
\end{center}

\medskip
\no $B(3,2)$-ben a következő 24 Hamilton-kör van:

\no 0010211220	\ 0020122110 \  0010221120 \ 0020112210 \  0011021220	\ 0022012110

\no 0011022120 \ 	0022011210 \ 0011202210 \ 0022101120 \  0011210220	\ 0022120110

\no 0011220210	\ 0022110120 \ 0011221020	\ 0022112010 \ 
0012022110	\ 0021011220

\no 0012110220	\ 0021220110 \ 0012202110	\ 0021101220 \ 
0012211020	\ 0021122010

\medskip\no A következő $B(3,2)$ gráfban két irányított Hamilton kört rajzoltunk be, az egyik élei zöldek, a másiké pirosak.

\includegraphics[scale=0.0625]{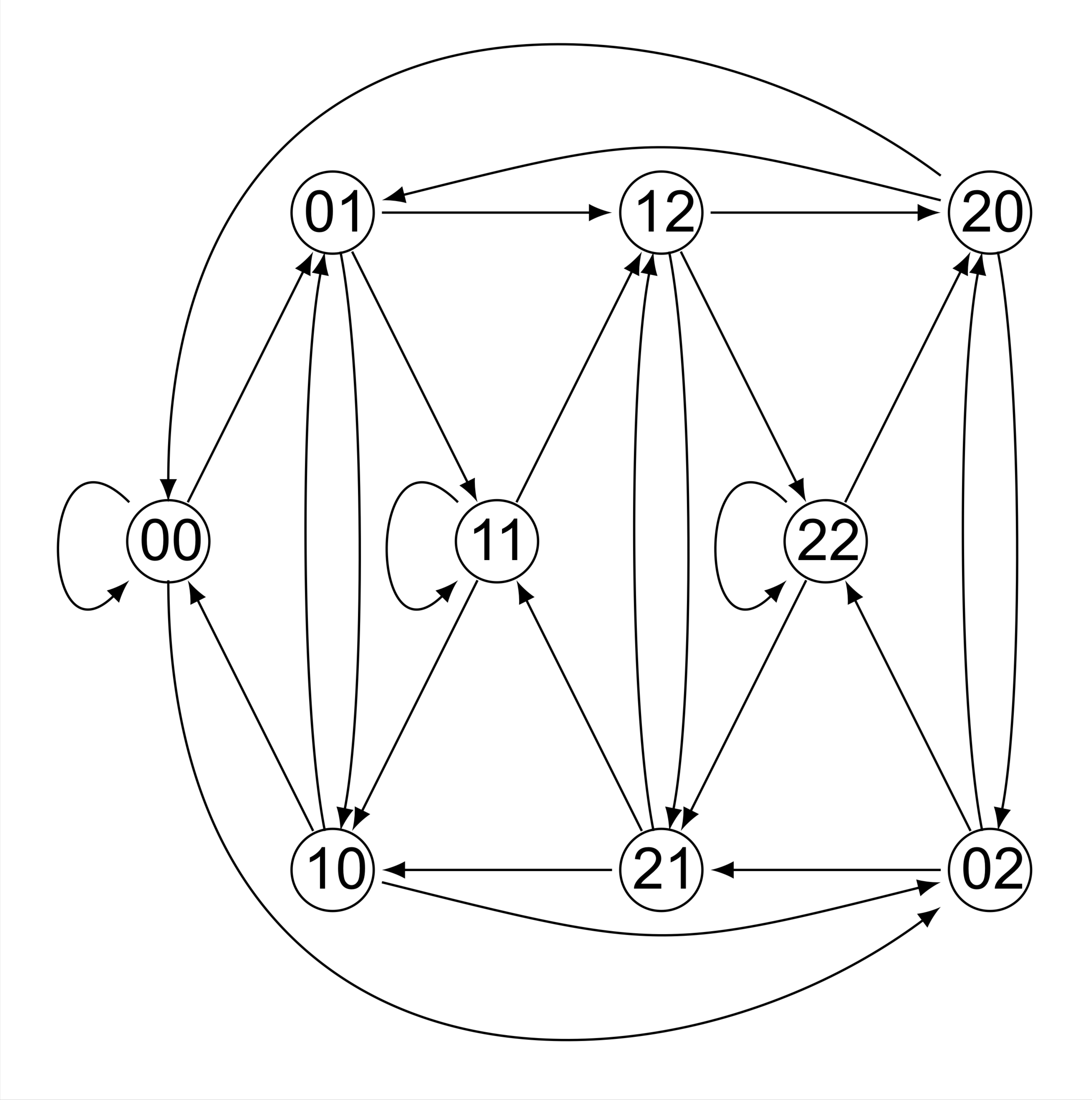} \  \includegraphics[scale=0.4]{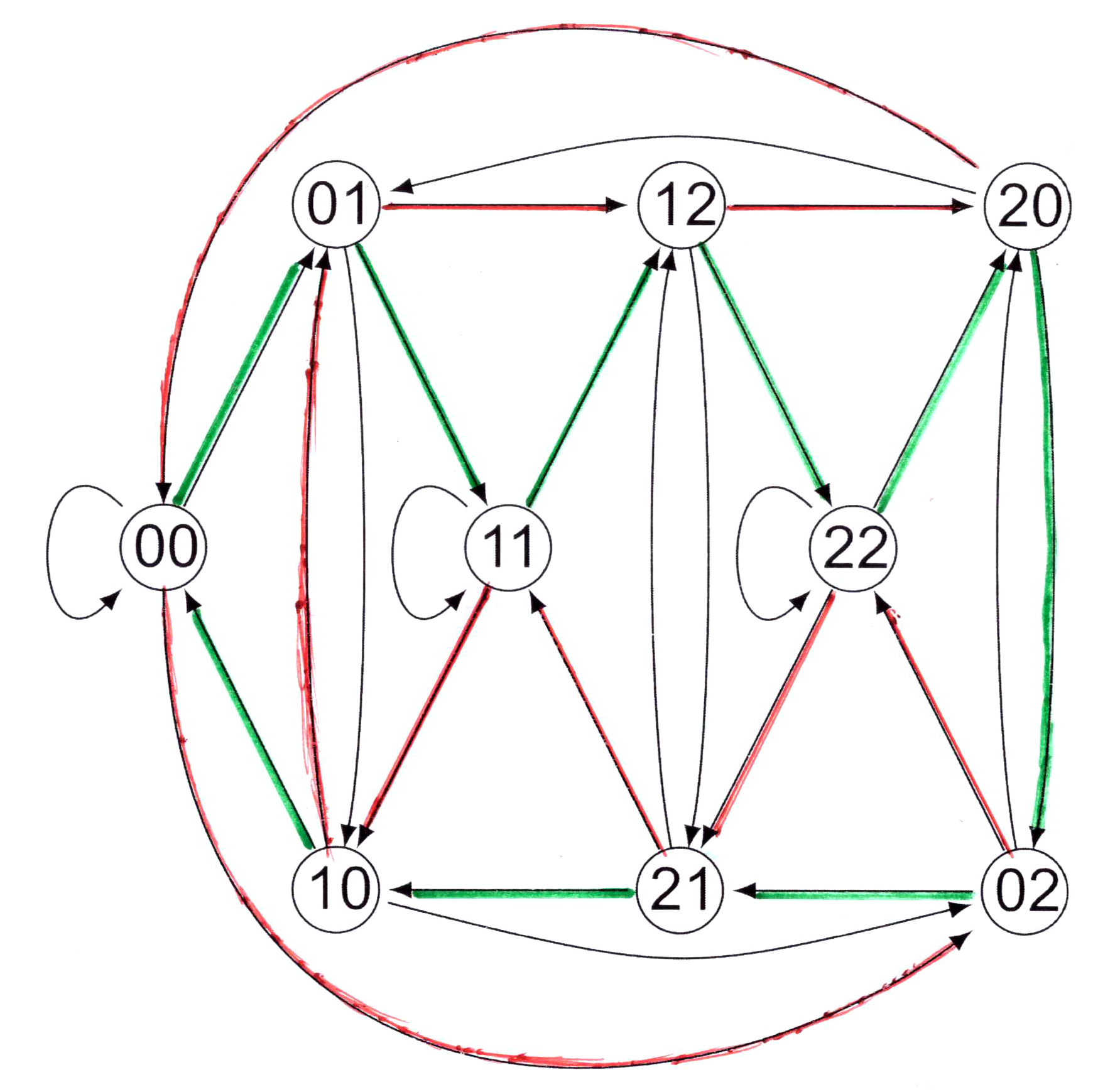}

\centerline{B(3,2)}

\medskip\no Láthatjuk, hogy a két irányított Hamilton-körnek nincsenek közös élei, az ilyen esetben azt mondjuk, hogy ez két Hamilton-kör élfüggetlen.

\no Mivel rengeteg Hamilton-kör van, feltevődik a kérdés, hogy ezek közül hány páronként élfüggetlen Hamilton-kör létezik. Látható, hogy a $B(3,2)$-ben két  élfüggetlen Hamilton-kör van. Több nem lehet, mert a $00$, $11$ és $22$ csúcsokból két él fut ki más csúcsokba, egy pedig önmagába, tehát ezeken csak két Hamilton-kör mehet át.   
Hasonló megfontolás alapján, a $B(n,m)$ gráfban is legfeljebb $n-1$ páronként élfüggetlen Hamilton-kör létezhet. Az a kérdés, hogy van-e ennyi? Johny Bond\footnote{Johny Bond romániai származású francia matematikus.} és Iványi Antal\footnote{Iványi Antal (1942--2017) magyar informatikus, egyetemi tanár az ELTE-n.} sejtése alapján pontosan ennyi van.

\medskip\no 
\keret{Johny Bond és Iványi Antal sejtése: \\
\kek{Ha $n\ge 2$ és $m\ge 2$, akkor a $B(n,m)$ gráfban  létezik $n-1$ páronként élfüggetlen irányított Hamilton-kör.}}

\no Egyelőre a sejtést nem sikerült bizonyítani, csak annyit, hogy legalább $\displaystyle\frac{n}{2}$ élfüggetlen Hamilton-kör létezik minden $B(n,m)$ gráfban.

\no Elmondhatjuk, hogy a rendezett De Bruijn-gráfokban a sok-sok rendezetlenül elhelyezkedő Hamilton-kör között létezik bizonyos számú (a sejtés szerint $n-1$) olyan irányított Hamilton-kör, hogy azok nem tartalmaznak közös éleket.

\no A $B(5,2)$ gráfban a következő Hamilton-körök élfüggetlenek:

\medskip
\kek{00102112041422430332313440}

\kek{00203223012133140443424110}

\kek{00304334023244210114131220}

\kek{00401441034311320221242330}

\smallskip
\no Meg lehet figyelni, hogy a második Hamilton-kört úgy kaptuk meg az elsőből, hogy 1 helyett 2-t, 2 helyett 3-t, 3 helyett 4-et, 4 helyett pedig 1-et írtunk. Hasonlóan képezzük a következőket is mindig az előzőből. Az  a kérdés, hogy az így kapott Hamilton-körök mindig élfüggetlen Hamilton kört adnak-e. Könnyű belátni, hogy ez nem igaz mindig, például $B(3,2)$ gráfban 
\kek{0010211220} és \kek{0020122110} mindegyike Hamilton-körnek felel meg, de van közös élük, például az 122-nek és a 211-nek megfelelő élek, azaz az 12 és 22, valamint a 21 és 11 csúcsok közöttiek. Az alábbi ábrán pirossal rajzoltuk az egyik, kékkel pedig a másik Hamilton-kör éleit. A szaggatott vonallal rajzoltak pirosak is, meg kékek is.

\begin{center}
 \begin{tikzpicture}[rotate=90, scale=3.3]
 \path (0:1cm)  node (X9){$00$};\draw (X9) circle  (3pt);
 \path (360/9*1:1cm)  node (X1){$01$};
 \draw (X1) circle  (3pt);
 \path (360/9*2:1cm)  node (X2){$10$};\draw (X2) circle  (3pt);
 \path (360/9*3:1cm)  node (X3){$02$};\draw (X3) circle  (3pt);
 \path (360/9*4:1cm)  node (X4){$21$};\draw (X4) circle  (3pt);
 \path (360/9*5:1cm)  node (X5){$11$};\draw (X5) circle  (3pt);
 \path (360/9*6:1cm)  node (X6){$12$};\draw (X6) circle  (3pt);
 \path (360/9*7:1cm)  node (X7){$22$};\draw (X7) circle  (3pt);
 \path (360/9*8:1cm)  node (X8){$20$};\draw (X8) circle  (3pt);
 
\draw[->, color=red, -latex] (X9) -- (X1);
\draw[->, color=red, -latex] (X1) -- (X2);
\draw[->, color=red, -latex] (X2) -- (X3);
\draw[->, color=red, -latex] (X3) -- (X4);
\draw[->, color=red, -latex] (X4) -- (X5);
\draw[->, color=red, -latex] (X5) -- (X6);
\draw[->, color=red, -latex] (X6) -- (X7);
\draw[->, color=red, -latex] (X7) -- (X8);
\draw[->, color=red, -latex] (X8) -- (X9);

\draw[->, color=blue, -latex] (X9) -- (X3);
\draw[->, color=blue, -latex] (X3) -- (X8);
\draw[->, color=blue, -latex] (X8) -- (X1);
\draw[->, color=blue, -latex] (X1) -- (X6);
\draw[->, color=blue, -latex, dashed] (X6) -- (X7);
\draw[->, color=blue, -latex] (X7) -- (X4);
\draw[->, color=blue, -latex, dashed] (X4) -- (X5);
\draw[->, color=blue, -latex] (X5) -- (X2);
\draw[->, color=blue, -latex] (X2) -- (X9);
 \end{tikzpicture}
\end{center}

\no Az viszont elképzelhető, hogy minden esetben vannak ilyen módon képzett független Hamilton-körök. Ha a $B(n,m)$ gráfban adott egy $H_1$ Hamilton-kör, képezzük a 
következő Hamilton-köröket:

\no $H_i$-t úgy kapjuk meg $H_{i-1}$-ből ($i=2,3,\ldots,n-1$), hogy elvégezzük  a következő cseréket:

$1 \rightarrow 2$
 
$2 \rightarrow 3$

$\cdots$

$n-2 \rightarrow n-1$

$n-1 \rightarrow 1$ 

és a $0$ nem változik.

\no Megfogalmaztam a következő sejtést.

\keret{
Sejtés

\kek{Bármilyen $n>2, m>1$ értékekre a $B(n,m)$ gráfban létezik egy olyan $H_1$ Hamilton-kör úgy, hogy a fenti eljárással kapott $H_2, H_3, \ldots , H_{n-1}$ Hamilton-körök a $H_1$-gyel együtt  páronként élfüggetlenek.}
}

\medskip\no Mivel a sejtést csak visszalépéses módszerrel tudom tesztelni, ez csak kis értékekre sikerült eddig, éspedig  a következőkre (néhol két megoldást is találtam):

\begin{small}
\bigskip\no
$\begin{array}{|c|c|l|}\hline
n & m& H_1 \\ \hline 
3 & 2 & 0011220210\\
3 & 2 & 0021011220\\
3 & 3 & 00010021011022202012111221200\\
4 & 2 & 00102113230331220\\
4 & 2 & 00102313033211220\\
4 & 3 & 000100210110201202310301311121130221232031323003332133122330322200\\
4 & 3 & 000100210110201202310301311121130223323003132123203330322213312200\\
5 & 2 & 00102112041422430332313440\\
6 & 2 & 0010211204131403325235505154534422430\\
7 & 2 & 00102112041306140315055162252353436442463326545660\\ \hline
\end{array}$  
\end{small}
	
\medskip\no Hamilton-kört a $B(n,m)$ gráfban könnyen elő lehet állítani az ún. Martin-algoritmussal, amely tulajdonképpen De Bruijn-szót generál. Elindulunk $m$ nullával. Legyen $s_1s_2\cdots s_k$ az eddig generált betűsorozat (a legelején $k=m$), ekkor a $(k+1)$-edik  pozícióban  folytatásként próbálkozunk $(n-1)$-nel (a betűk 0-tól $(n-1)$-ig vannak jelölve). Megvizsgáljuk, hogy ha az előző részben (azaz $s_1s_2\cdots s_k$-ban) nem szerepel már $s_{k-m+2}s_{k-m+3}\cdots s_ks_{k+1}$, akkor jó az $(n-1)$, de ha igen, akkor folytatjuk eggyel kisebb értékkel, és így tovább. Ha már nem tudjuk folytatni, az algoritmus befejeződik, és az eredmény a megfelelő De Bruijn-szó. Például: $n=3, m=2$-re ezt kapjuk:  \kek{$0022112010$}. (Pl. a 8. pozícióban nem lehet 2, mert előtte már van 22, de 1 sem lehet, mivel 21 is van, ezért kell 0. Az 11. pozícióba már egyik sem jó, mert a 02, 01 és 00 is szerepel már az előző részben. Ezért az algoritmus véget ér.)

\no Sajnos ez az algoritmus soha nem generál olyan $H_1$ Hamilton-kört, amelyre az előbbi módszerrel előállított $H_2, H_3, \ldots, H_{n-1}$ Hamilton-körök élfüggetlenek legyenek. De elképzelhető, hogy a Martin-algoritmushoz hasonló, de más algoritmus képes erre. Ez egyben bizonyítaná is a második sejtést, amelyet más módszerrel igen nehéz lenne bizonyítani.

\medskip\no 
A sejtést, hogy a $B(n,m)$ gráfban $n-1$ élfüggetlen Hamiltonion-kör van, a következő ábrák is sugalmazhatják. Mindegyik gráf úgy néz ki, mint egy virág, amelynek $n-1$ szirma van. 

\newpage
\begin{center}
{\qquad  De Bruijn-gráfok -- hurkok és irányítások elhagyásával}

\medskip
\noindent\includegraphics[scale=0.75]{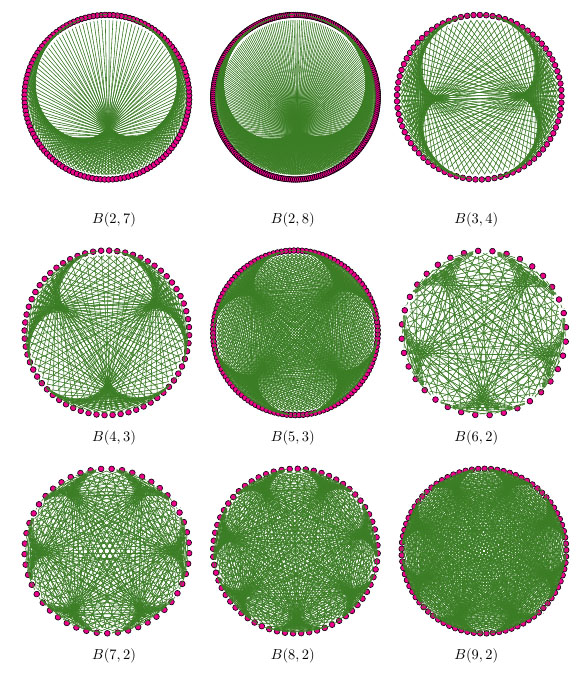} 

\end{center}

\centerline{\kek{\textbf{Szakirodalom}}}

\no Ha valaki érdeklődik az írásban foglalt fogalmak és eredmények felől, sokmindent megtalál a wikipédiában. Ha  a magyarban nem, akkor az angolban szinte biztosan.

\no A következő könyveket és cikkeket használtam.

\medskip
\no\textbf{1.} Andrásfai Béla: \textit{Ismerkedés a gráfelmélettel}, Tankönyvkiadó, Budapest, 1971.

\no\textbf{2.} Claude Berge: \textit{Théorie des graphes et ses applicatuions}, Dunod, Paris, 1967. (románul: \textit{Teoria grafurilor \cb si applica\cb tiile ei,} Editura Tehnică, Bucure\cb sti, 1969)

\no\textbf{3.} Johny Bond, Antal Iványi: Modelling of interconnection networks \mbox{using} De Bruijn graphs, in: \textit{Third Conference of Program Designer}, Ed. A. Iványi, Budapest, 1987, pp. 75--87.  \url{https://www.acta.sapientia.ro/PD/Third-conference-of-program-designers.pdf}

\no\textbf{4.} Zoltán Kása: \textit{On arc-disjoint Hamiltonian cycles in De Bruijn graphs}, \url{https://arxiv.org/abs/1003.1520}

\no\textbf{5.} Kása Zoltán, Anisiu Mira-Cristiana: Szavak bonyolultsága, in: Iványi Antal (szerk.): \textit{Informatikai algoritmusok III.,} Budapest, mondAt Kiadó, 2013, pp. 1697?1739., ISBN 9789638759689 (angolul: \url{https://www.researchgate.net/publication/274735246_Complexity_of_words}) 

\no\textbf{6.} Yaw-Ling Lin, Charles Ward, Bharat Jain,  Steven Skiena: Constructing Orthogonal de Bruijn Sequences, in \textit{Algorithms and Data Structures,} Lecture Notes in Computer Science, 2011, Volume 6844/2011, pp. 595--606.

\no\textbf{7.} Ioan Tomescu: \textit{Introducere în combinatorică,} Editura Tehnică, Bucure\cb sti, 1972 (magyarul: \textit{Kombinatorika és alkalmazása,} Budapest, 1978, fordította Maurer Gyula).

\medskip\no 2021. december 2.
     
 \newpage
 \centerline{\piros{\textbf{FÜGGELÉK}}}
 
 \bigskip
 \no A \kek{gráf} csúcsokból (szögpontokból, pontokból) és élekből áll. Egy él két csúcsot köt össze  egyenes szakasszal vagy görbe vonallal. A csúcsokat meg szoktuk jelölni, hogy könnyen tudjunk hivatkozni rájuk. Ha az éleknek van irányításuk, akkor \kek{irányított gráfokról} beszélünk. A bal oldali ábra egy (nem irányított) gráf (csúcsai 1, 2, 3, 4, 5), míg a jobb oldali irányított gráf (csúcsai $v_1, v_2, v_3, v_4$).   Az élek jelölése például: $\{3, 4\}$ vagy $(v_3,v_4)$ (itt számít  asorrend, ezért van egyszerű zárójel). 
 
 \bigskip\medskip
 \begin{center}
 \noindent\begin{tikzpicture}[scale=2]
 \tikzstyle{every node}=[draw,shape=circle];
 \path (0.5,2) node (v1) {\,1\,};
 \path (0,3) node (v2) {\,2\,};
 \path (1,3) node (v3) {\,3\,};
 \path (2,4) node (v4) {\,4\,};
 \path (2,3) node (v5) {\,5\,};
 \draw  
 (v1) -- (v2)
 (v2) -- (v3)     
 (v3) -- (v4)
 (v1) -- (v5)
 (v2) -- (v4)
 (v3) -- (v1)
 (v4) -- (v5) ;     
 \end{tikzpicture}
 \qquad\qquad
 \begin{tikzpicture}[scale=2]
 \tikzstyle{every node}=[draw,shape=circle];
 \path (0.5,3) node (v1) {\,$v_1$\,};
 \path (0,4) node (v2) {\,$v_2$\,};
 \path (1,3.7) node (v3) {\,$v_3$\,};
 \path (2,5) node (v4) {\,$v_4$\,};
 
 \draw [->, -latex] (v1) -- (v2);
 \draw [->,  -latex] (v2) -- (v3)     ;
 \draw [->,  -latex] (v3) -- (v4);
 \draw [->,  -latex] (v1) .. controls (1.75,3.5) ..  (v4);
 \draw [->,  -latex] (v2) -- (v4);
 \draw [->,  -latex] (v3) -- (v1) ;     	
 \end{tikzpicture}
 \end{center}
 	
\bigskip\medskip
\no A következő példákon illusztráljuk a séta, út és kör fogalmát. A \kek{séta} egymásutáni élek olyan sorozata, amelyben a csúcsok és élek is ismétlődhetnek is. Példa sétára: $\{v_6,v_7\}, \{v_7, v_3\}, \{v_3, v_1\}, \{v_1, v_2\}, \{v_2,v_3\}, \{v_3, v_4\}$, azaz röviden:  $v_6, v_7, v_3, v_1,v_2,v_3,v_4$ (az ábrán pirossal jelölve). Az \kek{út} olyan sajátos séta, amelyben nem ismétlődhetnek a csúcsok (és akkor nyilván az élek sem). Példa útra:  $v_6, v_1, v_5, v_4$ (az ábrán kékkel jelölve). A \kek{kör} olyan általánosított út, amelyben a kezdő- és a végcsúcsok azonosak. Példa:  $v_1, v_3, v_2, v_1$. Ezek a fogalmak kiterjeszthetők az irányított gráfokra is. Ekkor az élek csak az ugyanolyan irány szerint kapcsolódhatnak. Például a $(v_2,v_4)$, $(v_3,v_4)$ nem irányított út.

	\noindent\begin{tikzpicture}[scale=2.5]
	\tikzstyle{every node}=[draw,shape=circle];
	\path (0.5,2) node (v1) {\,$v_1$\,};
	\path (0,3) node (v2) {\,$v_2$\,};
	\path (1,3) node (v3) {\,$v_3$\,};
	\path (2,4) node (v4) {\,$v_4$\,};
	\path (2,3) node (v5) {\,$v_5$\,};
	\path (1,1.5) node (v6) {\,$v_6$\,};
	\path (2,1.5) node (v7) {\,$v_7$\,};
	\draw [color=red, thick] (v1) -- (v2) ;
	\draw [color=red, thick] (v2) -- (v3) ;
	\draw [color=red, thick]	(v3) -- (v4);
	\draw [color=blue, thick] (v1) -- (v5);
	\draw 	(v2) -- (v4);
	\draw [color=red, thick]	(v3) -- (v1);
	\draw [color=blue, thick]	(v4) -- (v5) ;
	\draw [color=blue, thick] (v6) -- (v1);
	\draw [color=red, thick]	(v6) -- (v7);
	\draw 	(v5) -- (v7);
	\draw [color=red]	(v3) -- (v7);     
	\end{tikzpicture} 
	\qquad \begin{tikzpicture}[scale=2.5]
\tikzstyle{every node}=[draw,shape=circle];
\path (0.5,2) node (v1) {\,$v_1$\,};
\path (0,3) node (v2) {\,$v_2$\,};
\path (1,3) node (v3) {\,$v_3$\,};
\path (2,4) node (v4) {\,$v_4$\,};
\path (2,3) node (v5) {\,$v_5$\,};
\path (1,1.5) node (v6) {\,$v_6$\,};
\path (2,1.5) node (v7) {\,$v_7$\,};
\draw [->, color=red, thick, -latex] (v2) -- (v1) ;
\draw [->, color=red, thick, -latex] (v3) -- (v2) ;
\draw [->, color=red, thick, -latex]	(v3) -- (v4);
\draw [->, color=blue, thick, -latex] (v1) -- (v5);
\draw [->, thick, -latex]	(v2) -- (v4);
\draw [->, color=red, thick, -latex]	(v1) -- (v3);
\draw [->, color=blue, thick, -latex]	(v5) -- (v4) ;
\draw [->, color=blue, thick, -latex] (v6) -- (v1);
\draw [->, color=red, thick, -latex]	(v6) -- (v7);
\draw [->, thick, -latex]	(v5) -- (v7);
\draw [->,color=red, thick, -latex]	(v7) -- (v3);     
\end{tikzpicture}

	\qquad $v_6, v_7, v_3, v_1,v_2,v_3,v_4$ \quad séta (illetve irányított séta)
	
	\medskip
	\qquad $v_6, v_1, v_5, v_4$ \quad  út (illetve irányított út)
	
	\medskip
	\qquad $v_1, v_3, v_2, v_1$ \quad  kör (illetve irányított kör)

\bigskip\bigskip

\no Ha egy $n$-csúcsú gráfban bármelyik két csúcs össze van kötve éllel, akkor azt a gráfot \kek{teljes gráfnak} nevezzük, és $K_n$-nel jelöljük. A $K_1$ egyetlen csúcsból áll, míg a $K_2$ két csúcsból és egy élből. Néhány példa teljes gráfokra: 

\centerline{\includegraphics[scale=0.35]{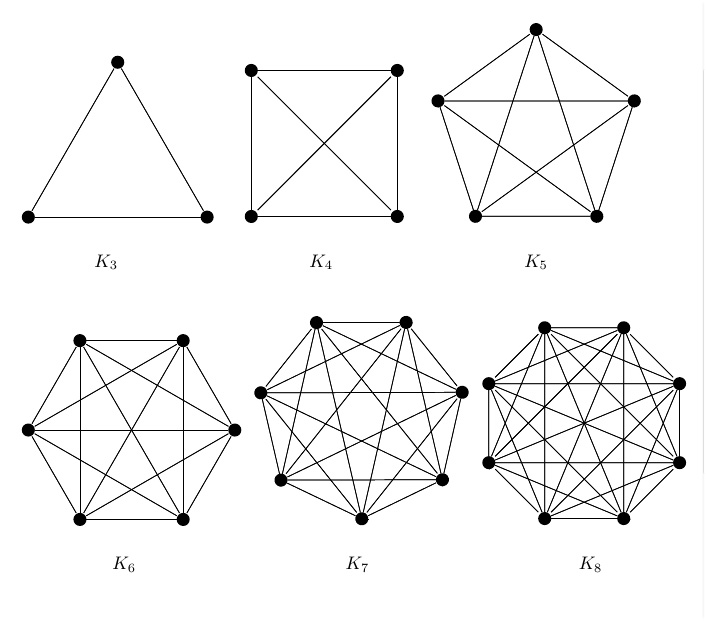}}	

A $G$ gráf \kek{komplementer gráfja} olyan gráf, amelyben a csúcsok ugyanazok mint $G$-ben, két csúcs között akkor és csak is  akkor van él, ha az eredeti gráfban nincs. Jelölés: $\overline{G}$.  A kettő együtt kiadja a teljes gráfot. A következő ábrán a komplementer gráfok piros élekkel vannak rajzolva. Egymásra helyezve, kiadják a teljes gráfot.

\bigskip\begin{tikzpicture}[scale=1]	
\pgfmathsetmacro{\n}{6};       
\foreach \i in {1,...,\n}
{  
	\path (360/\n*\i:1cm) node (X\i) { }; 
	\draw (X\i) circle  (3.5pt);               
}
\draw (X1) -- (X3);
\draw (X1) -- (X4);
\draw (X2) -- (X3);
\draw (X2) -- (X5);
\draw (X4) -- (X5);
\draw (X6) -- (X4);
\draw (X1) -- (X2);	
\draw (X6) -- (X5);
\end{tikzpicture}
\quad
\begin{tikzpicture}[scale=1]
\pgfmathsetmacro{\n}{6};       
\foreach \i in {1,...,\n}
{  
	\path (360/\n*\i:1cm) node (X\i) { }; 
	\draw (X\i) circle (3.5pt);               
}
\draw [color=red] (X1) -- (X5);
\draw [color=red] (X1) -- (X6);
\draw [color=red] (X2) -- (X4);
\draw [color=red] (X2) -- (X6);
\draw [color=red] (X3) -- (X4);
\draw [color=red] (X3) -- (X5);
\draw [color=red] (X3) -- (X6);
\end{tikzpicture}
\qquad\qquad
\begin{tikzpicture}[rotate=18,scale=1]
\pgfmathsetmacro{\n}{5};       
\foreach \i in {1,...,\n}
{  
	\path (360/\n*\i:1cm) node (X\i) { }; 
	\draw (X\i) circle (3.5pt);               
}

\draw (X1) -- (X2);
\draw (X2) -- (X3);
\draw (X3) -- (X4);
\draw (X4) -- (X5);
\draw (X5) -- (X1);	
\end{tikzpicture}
\quad
\begin{tikzpicture}[rotate=18,scale=1]
\pgfmathsetmacro{\n}{5};       
\foreach \i in {1,...,\n}
{  
	\path (360/\n*\i:1cm) node (X\i) { }; 
	\draw (X\i) circle (3.5pt);               
}

\draw [color=red] (X1) -- (X3);
\draw [color=red] (X1) -- (X4);
\draw [color=red] (X2) -- (X4);
\draw [color=red] (X2) -- (X5);
\draw [color=red] (X3) -- (X5);
\end{tikzpicture}

\bigskip 
\qquad\qquad $G$ \ \'es \ $\overline{G}$ \qquad\qquad\qquad\qquad \qquad\qquad\quad$G'$\ \'es \ $\overline{G'}$ 

\medskip
\begin{center}\begin{tikzpicture}[scale=1]	
\pgfmathsetmacro{\n}{6};       
\foreach \i in {1,...,\n}
{  
	\path (360/\n*\i:1cm) node (X\i) { }; 
	\draw (X\i) circle  (3.5pt);               
}
\draw (X1) -- (X3);
\draw (X1) -- (X4);
\draw (X2) -- (X3);
\draw (X2) -- (X5);
\draw (X4) -- (X5);
\draw (X6) -- (X4);
\draw (X1) -- (X2);	
\draw (X6) -- (X5);
\draw [color=red] (X1) -- (X5);
\draw [color=red] (X1) -- (X6);
\draw [color=red] (X2) -- (X4);
\draw [color=red] (X2) -- (X6);
\draw [color=red] (X3) -- (X4);
\draw [color=red] (X3) -- (X5);
\draw [color=red] (X3) -- (X6);	
\end{tikzpicture}
\qquad\qquad\qquad\qquad\qquad\qquad
\begin{tikzpicture}[rotate=18,scale=1]
\pgfmathsetmacro{\n}{5};       
\foreach \i in {1,...,\n}
{  
	\path (360/\n*\i:1cm) node (X\i) { }; 
	\draw (X\i) circle (3.5pt);               
}

\draw (X1) -- (X2);
\draw (X2) -- (X3);
\draw (X3) -- (X4);
\draw (X4) -- (X5);
\draw (X5) -- (X1);	

\draw[color=red] (X1) -- (X3);
\draw[color=red] (X1) -- (X4);
\draw[color=red] (X2) -- (X4);
\draw[color=red] (X2) -- (X5);
\draw[color=red] (X3) -- (X5);
\end{tikzpicture}
\end{center}

\no Elképzelhető olyan gráf is amelyben van hurokél (amelynek azonosak a  végpontjai) és  többszörös él (két csúcsot több él köt össze). Ha egy gráf nem ilyen, akkor \kek{egyszerű gráfnak} hívjuk.

\end{document}